\documentstyle[11pt,leqno]{article}
\newtheorem{theorem}{{\sc Theorem}}
\newcommand{\bt}{\begin{theorem}}
\newcommand{\et}{\end{theorem}}
\newcommand{\newsection}[1]{\setcounter{equation}{0} \setcounter{theorem}{0}
\section{#1}}

\newcommand{\NI}{\noindent}
\newcommand{\bea}{\begin{eqnarray}}
\newcommand{\eea}{\end{eqnarray}}

\def \spec#1 {\mathop{#1}}

\def \b #1 {\bf #1}

\newcommand {\CC}{\centerline}

\newcommand{\clf}{{\cal F}}

\newcommand{\ity}{\infty}
\newcommand{\raro}{\rightarrow}

\newcommand{\vsp}{\vskip 1em}
\newcommand{\vspp}{\vskip 2em}

\newcommand{\be}{\begin{equation}}
\newcommand{\ee}{\end{equation}}
\newcommand{\ben}{\begin{eqnarray*}}
\newcommand{\een}{\end{eqnarray*}}

\begin{document}

\sloppy
\CC {\bf  On the Long Range Dependence of Time-Changed}
\CC {\bf Generalized Mixed Fractional Brownian Motion }
\vsp
\CC {B.L.S. Prakasa Rao}
\CC{CR Rao Advanced Institute of Mathematics, Statistics}
\CC{and Computer Science, Hyderabad, India}

\vspp
\NI{\bf Abstract:} We introduce a generalized mixed fractional Brownian motion (gmfBm)  as a linear combination of two independent fractional Brownian motions with possibly different Hurst indices  and investigate conditions under which the time-changed gmfBm  exhibit long range dependence when the time-change is induced by a tempered stable subordinator or a Gamma process.
\vsp
\NI{\bf Keywords}: Fractional Brownian motion ; Mixed fractional Brownian motion; Generalized mixed fractional Brownian motion; Tempered  stable subordinator; Gamma process.
\vspp
\NI {\bf AMS Subject classification (2020)}: Primary 60G22.

\newsection{Introduction}

Geometric Brownian motion driven by a standard Brownian motion has been widely used for modeling fluctuations of share prices in a stock market. However efforts to model fluctuations in financial  markets with long range dependence through processes driven by a fractional Brownian motion were not successful as it was noted that such a model creates arbitrage opportunities contrary to the fundamental assumption of no arbitrage opportunity for model under rational market behaviour. Cheridito (2001) proposed modeling through processes driven by a mixed fractional Brownian motion. It was shown by Cheridito (2001) that a mixed fractional Brownian motion is a semimartingale if and only if the Hurst index $H$ is equal to $\frac{1}{2}$ reducing the process to a Wiener process or $H\in (3/4,1)$. Furthermore the probability measure generated by such a process is absolutely continuous with respect to the probability measure generated by a Wiener process if $H=1/2$ or $H\in (3/4,1).$ This in turn will lead to no arbitrage opportunities for modeling financial market behaviour through processes driven by a mixed fractional Brownian motion. This short discussion is just to motivate the study of processes driven by a mixed fractional Brownian motion.  

The problem of estimation of parameters for processes driven by processes which are mixtures of  independent Brownian and fractional Brownian motions started from the works of Cheridito (2001) and more recently in Prakasa Rao (2015a,b) among others. Mixed fractional Brownian models were  studied in Mishura (2008)  and Prakasa Rao (2010). A comprehensive review of fractional processes and their statistical inference was given in Prakasa Rao (2022). 

It is of interest to investigate sufficient conditions under which  some stochastic processes exhibit long range dependence for modeling stochastic phenomena such as  with the internet traffic, finance among other fields. Fractional Brownian motion (fBm) is one such process when the Hurst index of the process exceeds $\frac{1}{2}.$ It would be interesting if the class of stochastic processes with long range dependence is enlarged for stochastic modeling. Our aim in this paper is to introduce a new class of processes termed generalized mixed fractional Brownian motion (gmfBm) and study conditions under which the  time-changed gmfBm by either a tempered stable subordinator or a Gamma process has the long range dependence property. Kumar et al. (2019) investigated properties of a fBm delayed by  tempered and inverse tempered stable subordinator and Kumar et al. (2017) studied a fBm for long range dependence property when  it is time-changed by a Gamma process or an  inverse Gamma process. In a recent work, Alajmi and Mliki (2020,2021) studied the conditions under which a time-changed mixed fractional Brownian motion is long range dependent when the underlying process is a tempered stable subordinator or a Gamma process. Our  aim is to investigate conditions under which a gmfBm has the long-range dependence property when it is time-changed by a tempered stable subordinator or a Gamma process.
\vsp
Let $(\Omega, \clf, (\clf_t), P) $ be a stochastic basis satisfying the usual conditions. The natural filtration of a
stochastic process is understood as the $P$-completion of the filtration generated by this process. Let $B^{H}= \{B_t^H, t \geq 0 \}, $ be a  fractional Brownian motion (fBm) with the Hurst parameter $H \in (0,1)$, that is, a Gaussian process with continuous sample paths such that $B_0^H=0, E(B_t^H)=0$ and
\be
E(B_s^H B_t^H)= \frac{1}{2}[s^{2H}+t^{2H}-|s-t|^{2H}], t \geq 0, s \geq 0.
\ee
Let
$$N_t^{H_1,H_2}(a,b)= a B_t^{H_1}+ b B_t^{H_2}, t \geq 0.$$
where $B^{H_1}$ and $B^{H_2}$ are independents fractional Brownian motions and $a,b \in R$ not both equal to zero. The process $\{N_t^{H_1,H_2}(a,b), t \geq 0\}$ is called a {\it generalized mixed fractional Brownian motion} (gmfBm) with Hurst indices $H_1$ and $H_2.$
\vsp
The {\it time-changed generalized mixed fractional Brownian motion} is the process $Y_\beta^{H_1.H_2}(a,b)$ defined by
$$ Y_\beta^{H_1.H_2}(a,b) \equiv \{Y_{\beta_t}^{H_1.H_2}(a,b), t \geq 0\}=\{N_{\beta_t}^{H_1,H_2}(a,b), t \geq 0\}$$
where the process $N^{H_1,H_2}(a,b)$ is a gmfBm with parameters $a,b$ and $H_1,H_2 \in (0,1)$ and the process $\beta=\{\beta_t, t \geq 0\}$ is a subordinator independent of the fractional Brownian motions $B^{H_i}, i=1,2.$ A subordinator is a stationary process with independent increments.
\vsp
Let $s>0$ be fixed and $t>s.$ Suppose $\{X_t,t \geq 0\}$ is a stochastic process and  $ E(X_t^2)<\ity$ for all $t \geq 0.$ The process $\{X_t, t \geq 0\}$ is said to be {\it long-range dependent} if, for any fixed $s,$ 
$$ Corr(X_t,X_s) \simeq c(s)t^{-d}$$
as $t \raro \ity$ for some constant $c(s)$ depending only on $s.$
\vsp
We will now discuss the long range dependency properties of the time-changed process $Y_\beta^{H_1.H_2}(a,b)$ when the subordinator $\{\beta_t, t \geq 0\}$ is a tempered stable subordinator (TSS) or a Gamma process independent of the fractional Brownian motions $B^{H_i}, i=1,2.$.
\vsp
\newsection{GMFBM Time-changed by a Tempered Stable Subordinator}

A {\it tempered stable subordinator} (TSS) with index $\alpha \in (0,1)$ and tempering parameter $\lambda >0,$ is the non-decreasing non-negative Levy process $S^{\lambda,\alpha}= \{ S^{\lambda,\alpha}_t,t\geq 0\}$ where the random variable $S^{\lambda,\alpha}_t$ has the probability density function
$$f_{\lambda,\alpha}(x,t)=\exp(-\lambda x+\lambda^\alpha t)f_\alpha(x,t), \lambda >0, \alpha \in (0,1),$$
whe
$$f_\alpha(x,t)= \frac{1}{\pi}\int_0^\ity e^{-xy}e^{-ty^\alpha \cos \alpha \pi} \sin(t y^\alpha \sin \alpha \pi)dy.$$
\vsp
\NI{\bf Lemma 2.1:} For $q>0,$ 
$$E[(S_t^{\lambda,\alpha})^q] \simeq (\alpha \lambda^{\alpha-1}t)^q$$
as $t \raro \ity.$

For a proof of Lemma 2.1 and for more details about the process TSS, see Kumar et al. (2017). For convenience, we denote $Y_{S_t^{\lambda,\alpha}}^{H_1,H_2}(a,b)$ by $Y_{S_t^{\lambda,\alpha}}^{H_1,H_2}$ in the following computations.

We now investigate sufficient conditions under which the gmfBm $N^{H_1,H_2}$ which is time-changed by the TSS $S^{\lambda,\alpha}$ is long-range dependent. Observe that
\ben
\lefteqn{Cov(Y_{S_t^{\lambda,\alpha}}^{H_1,H_2}(a,b), Y_{S_s^{\lambda,\alpha}}^{H_1,H_2}(a,b))}\\\nonumber
&=& \frac{1}{2} E[(Y_{S_t^{\lambda,\alpha}}^{H_1,H_2})^2+(Y_{S_s^{\lambda,\alpha}}^{H_1,H_2})^2-(Y_{S_t^{\lambda,\alpha}}^{H_1,H_2}-Y_{S_s^{\lambda,\alpha}}^{H_1,H_2})^2]\\\nonumber
&=& \frac{1}{2} E[(N_{S_t^{\lambda,\alpha}}^{H_1,H_2})^2+(N_{S_s^{\lambda,\alpha}}^{H_1,H_2})^2-(N_{S_t^{\lambda,\alpha}}^{H_1,H_2}-N_{S_s^{\lambda,\alpha}}^{H_1,H_2})^2]\\\nonumber
&=& \frac{1}{2}E[(a B_{S_t^{\lambda,\alpha}}^{H_1}+b B_{S_t^{\lambda,\alpha}}^{H_2})^2+(a B_{S_s^{\lambda,\alpha}}^{H_1}+b B_{S_s^{\lambda,\alpha}}^{H_2})^2]\\\nonumber
&&\;\;\;\; -\frac{1}{2} E[(a(B_{S_t^{\lambda,\alpha}}^{H_1}-B_{S_s^{\lambda,\alpha}}^{H_1})+(B_{S_t^{\lambda,\alpha}}^{H_2}-B_{S_s^{\lambda,\alpha}}^{H_2}))^2].
\een
Since the fractional Brownian motions have stationary increments, it follows that
\ben
\lefteqn{Cov(Y_{S_t^{\lambda,\alpha}}^{H_1,H_2}(a,b), Y_{S_s^{\lambda,\alpha}}^{H_1,H_2}(a,b))}\\\nonumber
&=& \frac{1}{2}E[(a B_{S_t^{\lambda,\alpha}}^{H_1}+b B_{S_t^{\lambda,\alpha}}^{H_2})^2+(a B_{S_s^{\lambda,\alpha}}^{H_1}+b B_{S_s^{\lambda,\alpha}}^{H_2})^2\\\nonumber
&&\;\;\;-(a B_{S_{t-s}^{\lambda,\alpha}}^{H_1}+b B_{S_{t-s}^{\lambda,\alpha}}^{H_2})^2]\\\nonumber
&=& \frac{1}{2}E[(a B_{S_{t}^{\lambda,\alpha}}^{H_1})^2+(b B_{S_{t}^{\lambda,\alpha}}^{H_2})^2+2abB_{S_{t}^{\lambda,\alpha}}^{H_1}B_{S_{t}^{\lambda,\alpha}}^{H_2}]\\\nonumber
&&+\frac{1}{2}E[(a B_{S_{s}^{\lambda,\alpha}}^{H_1})^2+(b B_{S_{s}^{\lambda,\alpha}}^{H_2})^2+2abB_{S_{s}^{\lambda,\alpha}}^{H_1}B_{S_{s}^{\lambda,\alpha}}^{H_2}]\\\nonumber
&&-\frac{1}{2}E[(a B_{S_{t-s}^{\lambda,\alpha}}^{H_1})^2+(b B_{S_{t-s}^{\lambda,\alpha}}^{H_2})^2+2abB_{S_{t-s}^{\lambda,\alpha}}^{H_1}B_{S_{t-s}^{\lambda,\alpha}}^{H_2}].\\\nonumber
\een
By the independence of the fBms' $B^{H_1}$ and $B^{H_2}$ and their independence of the TSS $S^{\lambda,\alpha}$, we get that
\ben
\lefteqn{Cov(Y_{S_t^{\lambda,\alpha}}^{H_1,H_2}(a,b), Y_{S_s^{\lambda,\alpha}}^{H_1,H_2}(a,b))}\\\nonumber
&=& \frac{a^2}{2}[E(B_{S_t^{\lambda,\alpha}}^{H_1})^2+E(B_{S_s^{\lambda,\alpha}}^{H_1})^2-E(B_{S_{t-s}^{\lambda,\alpha}}^{H_1})^2]\\\nonumber
&&+\frac{b^2}{2}[E(B_{S_t^{\lambda,\alpha}}^{H_2})^2+E(B_{S_s^{\lambda,\alpha}}^{H_2})^2-E(B_{S_{t-s}^{\lambda,\alpha}}^{H_2})^2]\\\nonumber
\een
by observing that
\ben
E[B_{S_t^{\lambda,\alpha}}^{H_1}B_{S_t^{\lambda,\alpha}}^{H_2}]&=& E[E(B_{S_t^{\lambda,\alpha}}^{H_1}B_{S_t^{\lambda,\alpha}}^{H_2}|S_t^{\lambda,\alpha})]\\\nonumber
&=& \int E[B_z^{H_1}B_z^{H_2}]F_{S_t^{\lambda,\alpha}}(dz)\\\nonumber
&=& 0\\\nonumber
\een 
where $F_{S_t^{\lambda,\alpha}}(.)$ is the distribution function of the random variable $S_t^{\lambda,\alpha}.$
Hence
\ben
\lefteqn{Cov(Y_{S_t^{\lambda,\alpha}}^{H_1,H_2}(a,b), Y_{S_s^{\lambda,\alpha}}^{H_1,H_2}(a,b))}\\\nonumber
&=& \frac {a^2}{2}E[B^{H_1}(1)]^2 [ E(S_t^{\lambda,\alpha})^{2H_1}+ E(S_s^{\lambda,\alpha})^{2H_1}-E(S_[t-s]^{\lambda,\alpha})^{2H_1}]\\\nonumber
&&+\frac {b^2}{2}E[B^{H_2}(1)]^2 [E(S_t^{\lambda,\alpha})^{2H_2}+ E(S_s^{\lambda,\alpha})^{2H_2}-E(S_[t-s]^{\lambda,\alpha})^{2H_2}].\\\nonumber
\een
Fix $s.$ Applying Lemma 2.1, it follows that, for large $t >s,$
\ben
\lefteqn{Cov(Y_{S_t^{\lambda,\alpha}}^{H_1,H_2}(a,b), Y_{S_s^{\lambda,\alpha}}^{H_1,H_2}(a,b))}\\\nonumber
&\simeq & \frac{a^2}{2}(\alpha \lambda^{\alpha-1})^{2H_1}t^{2H_1}(2H_1\frac{s}{t}+E(S_{s,\alpha}^{\lambda,\alpha})^{2H_1}t^{-2H_1}
+O(t^{-2})\\\nonumber
&&\;\;\;\; + \frac{b^2}{2}(\alpha \lambda^{\alpha-1})^{2H_2}t^{2H_2}(2H_3\frac{s}{t}+E(S_{s,\alpha}^{\lambda,\alpha})^{2H_2}t^{-2H_2}
+O(t^{-2})\\\nonumber
&\simeq & a^2H_1s(\alpha \lambda^{\alpha-1})^{2H_1}t^{2H_1-1}+b^2H_2s(\alpha \lambda^{\alpha-1})^{2H_2}t^{2H_2-1}\\\nonumber
\een
by arguments similar to those in Alajmi and Milki (2021). Hence we have the following result.
\vsp
\NI{\bf Theorem 2.1:} Let $\{Y_{S_t^{\lambda,\alpha}}^{H_1,H_2}(a,b), t \geq 0\}$ be a gmfBm time-changed by the process $S^{\lambda \alpha}.$
Then, for fixed $s$ and large $t >s,$
 $$Cov(Y_{S_t^{\lambda,\alpha}}^{H_1,H_2}(a,b), Y_{S_s^{\lambda,\alpha}}^{H_1,H_2}(a,b))\simeq a^2H_1s(\alpha \lambda^{\alpha-1})^{2H_1}t^{2H_1-1}+b^2H_2s(\alpha \lambda^{\alpha-1})^{2H_2}t^{2H_2-1}.$$
\vsp
The next result deals with asymptotic behaviour of the second moment of the increments of the process $ \{Y_{S_t^{\lambda,\alpha}}^{H_1,H_2}(a,b), t \geq 0\}$ for fixed $s$ and large $t.$
\vsp
\NI{\bf Theorem 2.2:} Let $\{Y_{S_t^{\lambda,\alpha}}^{H_1,H_2}(a,b), t \geq 0\}$ be a gmfBm time-changed by the process $S^{\lambda \alpha}.$
Then, for fixed $s$ and large $t >s,$
\bea
\lefteqn{E[(Y_{S_t^{\lambda,\alpha}}^{H_1,H_2}(a,b)-Y_{S_t^{\lambda,\alpha}}^{H_1,H_2}(a,b))^2]}\\\nonumber
&\simeq & a^2H_1(\alpha \lambda^{\alpha-1})^{2H_1}t^{2H_1}-2a^2H_1(\alpha \lambda^{\alpha-1})^{2H_1}t^{2H_1-1}+a^2H_1(\alpha \lambda^{\alpha-1})^{2H_1}s^{2H_1}\\\nonumber
&&\;\;\;+b^2H_2(\alpha \lambda^{\alpha-1})^{2H_2}t^{2H_2}-2b^2H_2(\alpha \lambda^{\alpha-1})^{2H_2}t^{2H_2-1}+b^2H_2(\alpha \lambda^{\alpha-1})^{2H_2}s^{2H_2}.\\\nonumber
\eea
\NI{\bf Proof:} Let $s>0$ be fixed and let $t>s$ with $t \raro \ity.$ Then 
\ben 
\lefteqn{E[(Y_{S_t^{\lambda,\alpha}}^{H_1,H_2}(a,b)-Y_{S_t^{\lambda,\alpha}}^{H_1,H_2}(a,b))^2]}\\\nonumber
&=& E[Y_{S_t^{\lambda,\alpha}}^{H_1,H_2}(a,b)]^2 + E[Y_{S_s^{\lambda,\alpha}}^{H_1,H_2}(a,b)]^2 \\\nonumber
&&\;\;\;+E[Y_{S_t^{\lambda,\alpha}}^{H_1,H_2}(a,b)Y_{S_s^{\lambda,\alpha}}^{H_1,H_2}(a,b)]\\\nonumber
&\simeq & a^2H_1(\alpha \lambda^{\alpha-1})^{2H_1}t^{2H_1}-2a^2H_1(\alpha \lambda^{\alpha-1})^{2H_1}t^{2H_1-1}+a^2H_1(\alpha \lambda^{\alpha-1})^{2H_1}s^{2H_1}\\\nonumber
&&\;\;\;+b^2H_2(\alpha \lambda^{\alpha-1})^{2H_2}t^{2H_2}-2b^2H_2(\alpha \lambda^{\alpha-1})^{2H_2}t^{2H_2-1}+b^2H_2(\alpha \lambda^{\alpha-1})^{2H_2}s^{2H_2}\\\nonumber
\een
by observing that $ E[Y_{S_t^{\lambda,\alpha}}^{H_1,H_2}(a,b)]=0, t \geq 0$ and 
$$Cov[Y_{S_t^{\lambda,\alpha}}^{H_1,H_2}(a,b),Y_{S_s^{\lambda,\alpha}}^{H_1,H_2}(a,b)]=E[Y_{S_t^{\lambda,\alpha}}^{H_1,H_2}(a,b)Y_{S_s^{\lambda,\alpha}}^{H_1,H_2}(a,b)].$$ 
\vsp
We will now investigate conditions on the Hurst indices $H_1,H_2$ under which the process $\{Y_{S_t^{\lambda,\alpha}}^{H_1,H_2}(a,b), t \geq 0\}$ has long-range dependent behaviour.
\vsp
\NI{\bf Theorem 2.3:} Let $N^{H_1,H_2}(a,b)$ be a gmfBm with $0<H_1\leq H_2<1.$ Let $S^{\lambda ,\alpha}$ be a TSS process with $\lambda >0$ and $0<\alpha <1.$ Then the time-changed gmfBm by the process $S^{\lambda ,\alpha}$ has long-range dependent property if $2H_1-H_2 <1.$
\vsp
\NI{\bf Proof:} Let $s>0$ be fixed and $t>s$ with $t \raro \ity.$ Then, following Theorem 2.2, we get that
\ben
\lefteqn{Corr(Y_{S_t^{\lambda,\alpha}}^{H_1,H_2}(a,b), Y_{S_s^{\lambda,\alpha}})}\\\nonumber
& \simeq & \frac{a^2H_1s(\alpha \lambda^{\alpha-1})^{2H_1}t^{2H_1-1}+b^2H_2s(\alpha \lambda^{\alpha-1})^{2H_2}t^{2H_2-1}}{\sqrt{a^2H_1(\alpha \lambda^{\alpha-1})^{2H_1}t^{2H_1}+b^2H_2(\alpha \lambda^{\alpha-1})^{2H_2}t^{2H_2}}\sqrt{E[(Y_{S_t^{\lambda,\alpha}}^{H_1,H_2}(a,b))^2]}}\\\nonumber
&\simeq & \frac{c_1t^{2H_1-1}+c_2t^{2H_2-1}}{\sqrt{d_1t^{2H_1}+d_2t^{2H_2}}}\\\nonumber
&\simeq & \frac{c_3t^{2H_1-1}+c_4t^{2H_2-1}}{\sqrt{t^{2H_2}}}\;\;(\mbox{since}\; H_2 > H_1)\\\nonumber
&\simeq & c_5 t^{2H_1-H_2-1}+ c_6 t^{H_2-1} \\\nonumber
\een
as $s<t\raro \ity$ where $c_i, i=1,..,6$ are constants depending on $a,b,s,\alpha, \lambda, H_1,H_2$ but not $t.$ Since $0<H_2<1$ and $2H_1-H_2 <1$ by hypothesis, it follows that the last term tends to zero as $t \raro \ity.$ Hence the the time-changed gmfBm by the process $S^{\lambda ,\alpha}$ has long-range dependent property.
\vsp
\NI{\bf Remarks:} Theorem 2.3 extends the result of Alajmi and Milki (2021) for mfBm to gmfBm and it gives a sufficient condition for the long-range dependence property depending on the Hurst indices $H_1,H_2$ with $H_1<H_2.$
\vsp
\newsection{GMFBM Time-changed by a Gamma Process}
A Gamma process $\Gamma=\{\Gamma_t, t\geq 0\}$ is a stationary independent increment process with the Gamma distribution for $\Gamma_{t+s}-\Gamma_s$ with the probability density function given by
$$f(x,t)= \frac{1}{\Gamma(t/\nu)}x^{(t/\nu)-1}e^{-x}, x>0; f(x,t)=0, x \leq 0$$
where $\nu >0.$
\vsp
\NI{\bf Lemma 3.1:} For any $q>0,$
$$E[\Gamma_t^q] \simeq (t \nu^{-1})^q$$
as $t \raro \ity.$
\vsp
For a proof of Lemma 3.1, see Kumar et al. (2017). 

Let $N^{H_1,H_2}(a,b)=\{N_t^{H_1,H_2}(a,b), t\geq 0\}$ be a gmfBm and let $\Gamma= \{\Gamma_t,t\geq 0\}$ be a $\Gamma$ subordinator. Define
$$Y_{\Gamma_t}^{H_1,H_2}= N_{\Gamma_t}^{H_1,H_2}(a,b)= a B_{\Gamma_t}^{H_1}+b B_{\Gamma_t}^{H_2}, a,b \in R$$
not both zero. We now investigate sufficient conditions under which the gmfBm $N^{H_1,H_2}(a,b)$ which is time-changed by the $\Gamma$  process is long-range dependent. 

For convenience, we denote $Y_{\Gamma_t}^{H_1,H_2}(a,b)$ by $Y_{\Gamma_s}^{H_1,H_2}$ in the following computations. Assume that $0<H_1 < H_2<1.$ Fix $s>0.$ and let $t>s$ with $t \raro \ity.$ Suppose $0<H_1<H_2<1.$ Applying arguments similar to those used in Section 2, it is easy to prove that
\ben
\lefteqn{Cov(Y_{\Gamma_t}^{H_1,H_2},Y_{\Gamma_s}^{H_1,H_2})}\\\nonumber
&\simeq & \frac{2a^2H_1s}{\nu^{2H_1}}t^{2H_1-1}+\frac{2b^2H_2s}{\nu^{2H_2}}t^{2H_2-1} \\\nonumber
&\simeq & c_1 t^{2H_1-1}+c_2 t^{2H_2-1}\\\nonumber
&\simeq & c_3 t^{2H_1-H_2-1},\\\nonumber
\een
\be
Var(Y_{\Gamma_t}^{H_1,H_2}) \simeq c_4 t^{2H_1}+c_5t^{2H_2}
\ee
and
\ben
\lefteqn{E[(Y_{\Gamma_t}^{H_1,H_2}-Y_{\Gamma_s}^{H_1,H_2})^2]}\\\nonumber
&\simeq & \frac{2a^2H_1}{\nu^{2H_1}}t^{2H_1}-\frac{4a^2H_1s}{\nu^{2H_1}}t^{2H_1-1}+\frac{2a^2H_1}{\nu^{2H_1}}s^{2H_1}\\\nonumber
&&\;\;\;\;+\frac{2b^2H_2}{\nu^{2H_2}}t^{2H_2}-\frac{4b^2H_2s}{\nu^{2H_2}}t^{2H_2-1}+\frac{2b^2H_2}{\nu^{2H_2}}s^{2H_2}.\\\nonumber
\een
Furthermore, for fixed $s$ and $t>s$ with $t \raro \ity,$ it follows that
\ben 
\lefteqn{Corr(Y_{\Gamma_t}^{H_1,H_2},Y_{\Gamma_s}^{H_1,H_2})}\\\nonumber
&\simeq & \frac{|a|(2H_1)^{1/2}s}{t^H_1\nu^{H_1} \sqrt{E(Y_{\Gamma_s}^{H_1,H_2}})^2}t^{2H_1-1}\\\nonumber
&& \;\;\;\; + \frac{|b|(2H_2)^{1/2}s}{t^H_2\nu^{H_2} \sqrt{E(Y_{\Gamma_s}^{H_1,H_2}})^2}t^{2H_2-1}.\\\nonumber
\een
Note that
\ben
E(Y_{\Gamma_s}^{H_1,H_2})^2 \simeq \frac{2a^2H_1s^{2H_1-1}}{\nu^{2H_1}}+\frac{2b^2H_1s^{2H_2-1}}{\nu^{2H_2}}.\\\nonumber
\een
Combining the above estimates, it follows that
\ben 
\lefteqn{Corr(Y_{\Gamma_t}^{H_1,H_2},Y_{\Gamma_s}^{H_1,H_2})}\\\nonumber
&\simeq & c_6t^{2H_1-H_2-1}+c_7t^{H_2-1}
\een
for fixed $s$ and $t >s$ with $t \raro \ity.$ The constants $c_i, i=1,\dots,7$ depend on $a,b, H_1,H_2$ and $\nu.$ Observe that the last term tends to zero as $t \raro \ity$ if $2H_1-H_2<1$ since $0< H_2 <1.$ We have now the following result.
\vsp
\NI{\bf Theorem 3.1:} Let $N^{H_1,H_2}(a,b)$ be a gmfBm with $0<H_1 < H_2<1.$ Let $\Gamma$ be the Gamma process with parameter $\nu >0$. Then the time-changed gmfBm by the process $\Gamma$ has long-range dependent property if $2H_1-H_2 <1.$
\vsp
\NI {\bf Acknowledgment:} This work was supported under the scheme ``INSA Senior Scientist" of the Indian National Science Academy at the CR Rao Advanced Institute of Mathematics, Statistics and Computer Science, Hyderabad 500046, India.
\vsp
\NI {\bf References}
\begin{description}

\item Alajmi, S. and Milki, E. (2020) On the mixed fractional Brownian motion time changed by inverse $\alpha$-stable subordinator, {\it Applied Mathematical  Sciences}, {\bf 14}, 755-763.

\item Alajmi, S. and Milki, E. (2021) On the long range dependence of time-changed mixed fractional Brownian motion model, arXiv:2102.10180v1 [math.PR] 18 Feb 2021.\\

\item Cheridito, P. (2001) Mixed fractional Brownian motion, {\it Bernoulli}, {\bf 7}, 913-934.

\item Kumar, A., Wylomanska, A., Polozanski, R. and Sundar, S. (2017) Fractional Brownian motion time-changed by gamma and inverse gamma process, {\it Physica A: Statistical Mechanics and its Applications}, {\bf 468}, 648-667.\\

\item Kumar, A., Gajda, G., and Wylomanska, A. (2019) Fractional Brownian motion delayed by tempered and inverse tempered stable subordinators, {\it Methodol. Comput. Appl. Probab.}, {\bf 21}, 185-202.\\

\item Mishura, Y. (2008) {\it Stochastic Calculus for Fractional Brownian Motion and Related Processes}, Springer: Berlin.

\item Prakasa Rao, B.L.S. (2010) {\it Statistical Inference for Fractional Diffusion Processes}, London: Wiley.

\item Prakasa Rao, B.L.S. (2015a) Option pricing for processes driven by mixed fractional Brownian motion with superimposed jumps, {\it Probability in the Engineering and Information Sciences}, {\bf 29}, 589-596.

\item Prakasa Rao, B.L.S. (2015b) Pricing geometric Asian power options under mixed fractional Brownian motion environment, {\it Physica A}, {\bf 446}, 92-99.

\item Prakasa Rao, B.L.S. (2022) Fractional processes and their statistical inference: an overview, {\it Journal of the Indian
Institute of Science}, {\bf 102}, 1145-1175.

\end{description}
\end{document}